\documentclass[preprint,1p,12pt]{elsarticle}

\makeatletter
\def\ps@pprintTitle{
  \let\@oddhead\@empty
  \let\@evenhead\@empty
  \let\@evenfoot\@oddfoot
}
\makeatother

\usepackage{amssymb,amsmath,color}
\usepackage{graphicx,enumitem}

\usepackage{url}
\urldef{\mailsa}\path|{alfred.hofmann, ursula.barth, ingrid.haas, frank.holzwarth,|
\urldef{\mailsb}\path|anna.kramer, leonie.kunz, christine.reiss, nicole.sator,|
\urldef{\mailsc}\path|erika.siebert-cole, peter.strasser, lncs}@springer.com|

\newtheorem{theorem}{Theorem}
\newtheorem{corollary}[theorem]{Corollary}
\newtheorem{lemma}[theorem]{Lemma}
\newtheorem{proposition}[theorem]{Proposition}
\newdefinition{definition}[theorem]{Definition }
\newdefinition{remark}[theorem]{Remark }
\newdefinition{example}[theorem]{Example }
\newproof{proof}{Proof}

\begin{document}

\begin{frontmatter}

\title{On the {solvability} of bipolar max-product fuzzy relation equations with the standard negation}

\author
{M. Eugenia Cornejo, David Lobo, Jes\'us Medina}

\address
{Department of Mathematics,
 University of  C\'adiz. Spain\\
\texttt{\{mariaeugenia.cornejo,david.lobo,jesus.medina\}@uca.es}}

\begin{abstract}
Bipolar fuzzy relation equations arise when unknown variables together with their logical negations appear simultaneously in fuzzy relation equations. This paper gives a characterization of the solvability of bipolar max product fuzzy (relation) equations with the standard negation. In addition, some properties associated with the existence of the greatest/least solution or maximal/minimal solutions are shown, when these (relation) equations are solvable. Different examples are included in order to clarify the developed theory.
\end{abstract}

\end{frontmatter}

\section{Introduction.}

{Situations and ideas with positive and negative effects are daily present in any human being's life. The principle behind any decision is a cautious examination of pros and cons, that is, of positive and negative sides of information. This inherent nature of human reasoning makes bipolarity a desirable feature to be represented in intelligent technologies.

This manuscript deals with a special kind of bipolarity, in which positive and negative aspects of information are independently evaluated on a single set of plausible values. In particular, negative aspects are precisely determined from positive aspects by means of a negation operator. According to the three forms of bipolarity presented in~\cite{duboisprade2008}, this work follows the philosophy of type II bipolarity. Hence, our notion of bipolarity is somehow associated with other frameworks in which type II bipolarity arises. Namely, type II bipolarity can be found in uncertainty theories dealing with incomplete information, such as belief and plausibility functions of Shafer~\cite{Shafer1976}, possibility and necessity measures~\cite{Dubois1988}, and intuitionistic fuzzy sets~\cite{Atanassov1999}.

Bipolar fuzzy relation equations were introduced in~\cite{Freson2013} as a generalization of fuzzy relation equations (FREs) in which unknown variables appear together with their logical negations simultaneously. Since its first definition in the 1980s by Elie Sanchez~\cite{sanchez76,sanchez79}, many papers have been published on the resolution of FREs~\cite{baets99_eq,dm:mare,Medina2016,Ignjatovic20151,Li2008,Medina2017:ija,Peeva2016,perfilieva08_fss}. In particular, essential results for solving max-product fuzzy relation equations can be found in~\cite{BOURKE1998,Loetamonphong:99,Peeva2013,Peeva2007}.}

As far as we know, the literature devoted to the resolution of bipolar fuzzy relation equations is really limited.
Bipolar fuzzy relation equations based on max-min composition are studied in~\cite{Li2016, Liu2015} and those based on max-product composition are analysed in~\cite{fuzzieee2017,escim2017,CLM:JCAM2018}. There are also a few papers which deal with the applicative perspective of bipolar equations and all of them are related to optimization problems~\cite{Freson2013,Li2014,Zhou2016}. Hence, a deeper study of these equations from different perspectives will be useful in order to increase their flexibility and applicability in real cases.

This paper provides new advances with respect to the initial study about bipolar max-product fuzzy equations with the standard negation, {firstly presented by the authors in~\cite{fuzzieee2017}.}
Hence, two of the most {common and practical} operators have been considered: the product t-norm, which  is also known as the Goguen t-norm; and the standard negation, which coincides with the residuated negation of the \L ukasiewicz t-norm.
{Notice that, the product t-norm has particularly interesting properties, giving rise to successful applications in the literature, for instance,~\cite{figueroa2018,gerami2019}.} On the other hand, the standard negation and involutive negations in general, have a great interest in the theory~\cite{Cintula2010,ciungu14,cemr:15} and in the applications~\cite{deer2018,dinola2015,Kochergin2019}.

We will start {this paper} in Section~\ref{sec:bipolarFE} including the definition of bipolar max-product fuzzy equation with the standard negation containing different unknown variables. The technical results related to the existence of the greatest/least solution or a finite number of maximal/minimal solutions for solvable bipolar equations will also be provided. It is important to mention that this study will also consider bipolar equations whose independent term is equal to zero, complementing the results given in~\cite{fuzzieee2017}. Then, in Section~\ref{sec:bipolarFRE}, a characterization for the solvability of an arbitrary bipolar max-product fuzzy relation equation with the standard negation will be introduced. {Different properties related to the existence and the analytical description of the maximal and minimal solutions of an arbitrary bipolar max-product fuzzy relation equation are given {in Section~\ref{sec:allsol}}}. Finally, this paper finishes in Section~\ref{sec:conclusion} with some conclusions and prospects for future work.  

\section{
Bipolar max-product fuzzy equations with the standard negation}\label{sec:bipolarFE}
This section carries out {a study} on the solvability of bipolar fuzzy   equations based on the max-product t-norm composition and the standard negation. Specifically, the standard negation $n\colon[0,1]\to[0,1]$ is defined as $n(x)=1-x$, for each $x\in[0,1]$. First of all, we will show the sufficient and necessary conditions which guarantee that bipolar max-product fuzzy equations containing different unknown variables are solvable. The following step of our study will consist in knowing when solvable bipolar fuzzy equations have a greatest (least, respectively) solution or a finite number of maximal (minimal, respectively) solutions.

Next, we will introduce formally the definition of bipolar max-product fuzzy  equation with the standard negation, and the characterization theorem {on its solvability}. {From now on, let us fix $m\in \mathbb N$ and $a_j^+,a_j^-,b\in [0,1]$, with $j\in \{1,\dots,m\}$.}

\begin{definition}
{Let $x_j\in[0,1]$ be an unknown value, for each $j\in \{1,\dots,m\}$, $\ast$   the product t-norm and $\vee$   the maximum operator. Equation~\eqref{eq:bipolarFRE} is called \emph{bipolar max-product fuzzy equation with the standard negation}.
\begin{equation}\label{eq:bipolarFRE}
\bigvee_{j=1}^m (a_j^+* x_j) \vee (a_j^-* (1-x_j))=b
\end{equation}
The \emph{corresponding max-product  fuzzy  equation} of Equation~\eqref{eq:bipolarFRE} is defined as
  \begin{equation}\label{eq:FRE}
    \bigvee_{j=1}^m (a_j^+* x_j) \vee (a_j^-* y_j)=b
  \end{equation}
where $x_j,y_j\in[0,1]$ are unknown values, for each $j\in \{1,\dots,m\}$.
}
\end{definition}

{Notice} that, if $a_j^+=a_j^-=0$ for a given $j\in\{1,\dots,m\}$, then we can remove the unknown value $x_j$ from the equation, since it may take any value in $[0,1]$. Hence, hereinafter, we assume that either $a_j^+\neq0$ or $a_j^-\neq0$ for each $j\in\{1,\dots,m\}$.

The following result {was advanced in~\cite{fuzzieee2017}, without proof, and} presents a characterization of the solvability of these equations. Specifically, it establishes that Equation~\eqref{eq:bipolarFRE} is solvable if and only if the components of the greatest solution of its corresponding max-product  fuzzy  equation verifies a certain inequality.

\begin{theorem}\label{th:luka}
{The bipolar max-product fuzzy equation~\eqref{eq:bipolarFRE} is solvable if and only if its corresponding max-product fuzzy equation~\eqref{eq:FRE} is solvable and the inequality $1\leq \bar{x}_j+ \bar{y}_j$ holds, for all $j\in\{1,\dots,m\}$, where $(\bar{x}_1,\bar{y}_1,\dots,\bar{x}_m,\bar{y}_m)\in[0,1]^{2m}$ is its greatest solution.}
\end{theorem}
\begin{proof}
Supposing that $1\leq \bar{x}_j+ \bar{y}_j$, for each $j\in\{1,\dots,m\}$, we will prove that Equation~\eqref{eq:bipolarFRE} is solvable. In order to reach this conclusion, consider the tuple $(\hat{x}_1,\dots,\hat{x}_m)$ defined as:
  \[\hat{x}_j= \left\{\begin{array}{ll}
             1-\bar{y}_j & \hbox{ if}\quad  a_j^-* \bar{y}_j=b\\
             \bar{x}_j & \hbox{ otherwise}
             \end{array}
  \right.\]
for each $j\in\{1,\dots,m\}$. We will see that $(\hat{x}_1,\dots,\hat{x}_m)$ is a solution of Equation~\eqref{eq:bipolarFRE}. 

Clearly, the tuple $(\bar{x}_1,\bar{y}_1,\dots,\bar{x}_m,\bar{y}_m)$ being a solution of Equation~\eqref{eq:FRE} implies that either there exists $k\in\{1,\dots,m\}$ such that $a_{k}^+* \bar{x}_{k}=b$, or there exists ${k}\in\{1,\dots,m\}$ such that $a_{k}^-* \bar{y}_{k}=b$ (clearly, both cases are also possible at the same time).  On the one hand, if there exists $k\in\{1,\dots,m\}$ such that $a_{k}^-* \bar{y}_{k}=b$, then by definition $\hat{x}_{k}=1-\bar{y}_{k}$. Observe that, as $(\bar{x}_1,\bar{y}_1,\dots,\bar{x}_m,\bar{y}_m)$ is a solution of Equation~\eqref{eq:FRE}, the inequality $a_{k}^+*\bar{x}_{k}\leq b$ is verified. Since $1\leq \bar{x}_{k}+ \bar{y}_{k}$ by hypothesis, or equivalently $1- \bar{y}_{k}\leq\bar{x}_{k}$, we obtain that the following chain of inequalities holds:
\[a_{k}^+* \hat{x}_{k}=a_{k}^+* (1-\bar{y}_{k})\leq a_{k}^+*\bar{x}_{k}\leq b\]
Consequently, we have that:
\[(a_{k}^+* \hat{x}_{k})\vee (a_{k}^-* (1-\hat{x}_{k}))=b\]
On the other hand, suppose that there exists $k\in\{1,\dots,m\}$ such that $a_{k}^+* \bar{x}_{k}=b$ and $a_{k}^-* \bar{y}_{k}\neq b$. Notice that, if $a_{k}^-* \bar{y}_{k}= b$, then we are in the previous case. By definition, $\hat{x}_{k}=\bar{x}_{k}$, and due to $(\bar{x}_1,\bar{y}_1,\dots,\bar{x}_m,\bar{y}_m)$ is a solution of Equation~\eqref{eq:FRE}, we obtain that {$a_{k}^-* \bar{y}_{k}\leq b$}. As a consequence, from the inequality $1\leq \bar{x}_{k}+ \bar{y}_{k}$, or equivalently $1- \bar{x}_{k}\leq\bar{y}_{k}$, we can assert that: 
\[a_{k}^-* (1-\hat{x}_{k})=a_{k}^-* (1-\bar{x}_{k})\leq a_{k}^-*\bar{y}_{k}\leq b\]
As a result, we obtain that:
\[(a_{k}^+* \hat{x}_{k})\vee (a_{k}^-* (1-\hat{x}_{k}))=b\]
In order to demonstrate that $(\hat{x}_1,\dots,\hat{x}_m)$ is a solution of Equation~\eqref{eq:bipolarFRE}, it remains to see that the next inequality holds, for each $j\in\{1,\dots,m\}$ with $j\neq k$:
\[(a_j^+* \hat{x}_j)\vee (a_j^-* (1-\hat{x}_j))\leq b\]
Given $j\in\{1,\dots,m\}$ with $j\neq k$, we will consider two cases: if $a_j^-* \bar{y}_j=b$, then following an analogous reasoning to the previous one $\hat{x}_{j}=1-\bar{y}_{j}$, and hence:
\[(a_{j}^+* \hat{x}_{j})\vee (a_{j}^-* (1-\hat{x}_{j}))=b\]
On the contrary, if $a_j^-* \bar{y}_j\neq b$, then by definition $\hat{x}_j=\bar{x}_j$. Since $(\bar{x}_1,\bar{y}_1,\dots,\bar{x}_m,\bar{y}_m)$ is a solution of Equation~\eqref{eq:FRE}, the inequalities $a_{j}^+*\bar{x}_j\leq b$ and $a_j^-*\bar{y}_j< b$ hold. As a result, taking into account that $1\leq \bar{x}_{j}+ \bar{y}_{j}$ by hypothesis, we can conclude that the following chain is satisfied:
\[(a_{j}^+* \hat{x}_j)\vee (a_{j}^-* (1-\hat{x}_j))=(a_{j}^+* \bar{x}_j)\vee (a_{j}^-* (1-\bar{x}_j))\leq (a_{j}^+* \bar{x}_j)\vee (a_{j}^-* \bar{y}_j)\leq b\]
Therefore, $(\hat{x}_1,\dots,\hat{x}_m)$ is a solution of Equation~\eqref{eq:bipolarFRE}.

In order to prove the other implication, we will assume that there exists $j\in\{1,\dots,m\}$ such that $\bar{x}_j+ \bar{y}_j<1$ and let us see that Equation~\eqref{eq:bipolarFRE} is not solvable. For that, we will suppose that there exists a solution $(\hat{x}_1,\dots,\hat{x}_m)$ of Equation~\eqref{eq:bipolarFRE}, and we will obtain a contradiction. 

{Clearly, as $(\hat{x}_1,\dots,\hat{x}_m)$ is a solution of Equation~\eqref{eq:bipolarFRE}}, then the tuple $(\hat{x}_1,1-\hat{x}_1,\dots,\hat{x}_m,1-\hat{x}_m)$ forms a solution of Equation~\eqref{eq:FRE}. Moreover, since the tuple $(\bar{x}_1,\bar{y}_1,\dots,\bar{x}_m,\bar{y}_m)$ is the greatest solution of Equation~\eqref{eq:FRE}, we obtain that, in particular, $\hat{x}_j\leq\bar{x}_j$ and $1-\hat{x}_j\leq \bar{y}_j$. Therefore, the following chain of inequalities holds $1\leq\hat{x}_j+ \bar{y}_j\leq\bar{x}_j+\bar{y}_j$. This fact leads us to a contradiction, due to $\bar{x}_j+ \bar{y}_j<1$ by hypothesis. Consequently, we can conclude that Equation~\eqref{eq:bipolarFRE} is solvable if and only if the inequality $1\leq \bar{x}_j+ \bar{y}_j$ holds, for all $j\in\{1,\dots,m\}$.
\qed\end{proof}

Notice that, if there exists $k\in\{1,\dots,m\}$ satisfying the equalities $a_k^+* \bar{x}_k= b$ and $a_k^-* \bar{y}_k=b$, then the tuple $(\hat{x}_1,\dots,\hat{x}_{k-1},(1-\bar{y}_k),\hat{x}_{k+1},\hat{x}_m)$ is also a solution of Equation~\eqref{eq:bipolarFRE}. In addition, we know that this solution is smaller than the solution obtained in the proof of Theorem~\ref{th:luka} since, by hypothesis, $(1-\bar{y}_i)\leq \bar{x}_i$, for each $i\in\{1,\dots,m\}$.

Although the algebraic structure corresponding to the complete solution set of Equation~\eqref{eq:bipolarFRE} has been studied in~\cite{fuzzieee2017} as well, the case in which the independent term is equal to zero was not considered. It is important to highlight that, if $b=0$, then the solution set of Equation~\eqref{eq:bipolarFRE} has a simple structure. Indeed, in that case, the solvability of Equation~\eqref{eq:bipolarFRE} can be characterized by means of the values of the coefficients $a_j^+$ and $a_j^-$, as the next proposition shows.

\begin{proposition}\label{prop:b=0}
{The bipolar max-product fuzzy equation
  \begin{equation}\label{eq:b=0}
    \bigvee_{j=1}^m (a_j^+* x_j) \vee (a_j^-* (1-x_j))=0
  \end{equation}
is solvable if and only if, for each $j\in\{1,\dots,m\}$, either $a_j^+=0$ or $a_j^-=0$.}
\end{proposition}
\begin{proof}
On the one hand, suppose that, for each $j\in\{1,\dots,m\}$, either $a_j^+=0$ or $a_j^-=0$. It is important to highlight that $a_j^+$ and $a_j^-$ cannot be zero simultaneously, since in such case the unknown value $x_j$ can be removed from the equation. Clearly, the tuple 
$(\hat{x}_1,\dots,\hat{x}_m)$ given by
\[\hat{x}_j= \left\{\begin{array}{ll}
             1 & \hbox{ if}\quad  a_j^+=0\\
             0 & \hbox{ if}\quad  a_j^-=0
             \end{array}
  \right.\]
is straightforwardly a solution of Equation~\eqref{eq:b=0}.

On the other hand, assume that Equation~\eqref{eq:b=0} is solvable and the tuple $(\hat{x}_1,\dots,\hat{x}_m)$ is a solution of Equation~\eqref{eq:b=0}. We will suppose that there exists $k\in\{1,\dots,m\}$ such that $a_k^+>0$ and $a_k^->0$ and we will obtain a contradiction. We will distinguish two cases:
\begin{itemize}
\item  If $x_k=0$, then the chain of inequalities below holds:
  \[\bigvee_{j=1}^m(a_j^+ \ast x_j) \vee (a_j^-\ast (1-x_j))\geq (a_k^-\ast (1-x_k))=a_k^->0\]
  which contradicts that  $(\hat{x}_1,\dots,\hat{x}_m)$ is  a solution of Equation~\eqref{eq:b=0}.
\item If $x_k>0$, we obtain the same strict inequality:
  \[\bigvee_{j=1}^m(a_j^+ \ast x_j) \vee (a_j^-\ast (1-x_j))\geq(a_k^+\ast x_k)>0\]
which also contradicts $(\hat{x}_1,\dots,\hat{x}_m)$ is  a solution of Equation~\eqref{eq:b=0}.
\end{itemize}
Thus,  we can ensure that either the equality $a_j^+=0$ or $a_j^-=0$ holds, for all $j\in\{1,\dots,m\}$.
\qed\end{proof}

Proposition~\ref{prop:b=0} leads us to assert that if Equation~\eqref{eq:b=0} is solvable then there is only one solution. This result is formalized in the next proposition.

\begin{proposition}\label{prop2:b=0}
{If Equation~\eqref{eq:b=0} is solvable, then the tuple $(\hat{x}_1,\dots,\hat{x}_m)$ defined as:
\[
\hat{x}_j= \left\{\begin{array}{ll}
1 & \hbox{ if}\quad a_j^+=0\\
             0 & \hbox{ if}\quad  a_j^-=0
\end{array}
\right.
\]
for each $j\in\{1,\dots,m\}$, is the unique solution of Equation~\eqref{eq:b=0}.}
\end{proposition}
\begin{proof}
{First of all, taking into account Proposition~\ref{prop:b=0}, we can assert that the tuple $(\hat{x}_1,\dots,\hat{x}_m)$ defined as
\[\hat{x}_j= \left\{\begin{array}{ll}
1 & \hbox{ if}\quad a_j^+=0\\
             0 & \hbox{ if}\quad  a_j^-=0
\end{array}
\right.\]
for each $j\in\{1,\dots,m\}$ is well-defined. Additionally, following an analogous reasoning to the proof given {to} Proposition~\ref{prop:b=0}, we obtain that $(\hat{x}_1,\dots,\hat{x}_m)$ is a solution of Equation~\eqref{eq:b=0}.} 

Now, suppose that it is not its unique solution. As a consequence, there exists a solution $(x_1,\dots,x_m)$ of Equation~\eqref{eq:b=0} such that $(x_1,\dots,x_m)\neq(\hat{x}_1,\dots,\hat{x}_m)$. Thus, one of the following statements holds:
\begin{itemize}
\item[(a)] there exists $k\in\{1,\dots,m\}$ satisfying $\hat{x}_k<x_k$.
\item[(b)] there exists $k\in\{1,\dots,m\}$ satisfying $x_k<\hat{x}_k$.
\end{itemize}
On the one hand, suppose that (a) is verified. If $a_k^+=0$  then $\hat{x}_k=1$, by definition. Since $x_k\in[0,1]$, the inequality $x_k\leq \hat{x}_k$ holds, which contradicts (a). As a consequence, $a_k^+>0$ and this fact implies that $0=\hat{x}_k< x_k$. These two inequalities lead us to the following chain:
\[
0< a_k^+* {x}_k \leq (a_k^+* {x}_k) \vee (a_k^-* (1-{x}_k)) 
\]
Hence, $(x_1,\dots,x_m)$ is not a solution of Equation~\eqref{eq:b=0} in contradiction with the hypothesis. On the other hand, if (b) holds, a dual process leads us to a contradiction. Therefore, we conclude that $(\hat{x}_1,\dots,\hat{x}_m)$ is the unique solution of Equation~\eqref{eq:b=0}.
\qed\end{proof}

In the following, we will show under what conditions a bipolar max-product fuzzy  equation, defined with the standard negation and whose independent term is different from zero, has either a greatest (least, respectively) solution or a finite number of maximal (minimal, respectively) solutions. {This result was previously advanced in~\cite{fuzzieee2017} without proof.}

\begin{theorem}\label{th:luka2}
{Let Equation~\eqref{eq:bipolarFRE} be solvable, $b\neq0$ and $(\bar{x}_1,\bar{y}_1,\dots,\bar{x}_m,\bar{y}_m)\in[0,1]^{2m}$ be the greatest solution of its corresponding max-product fuzzy equation~\eqref{eq:FRE}. The following statements hold:
\begin{itemize}
\item[(a)] If there exists some $k\in\{1,\dots,m\}$ such that $a_k^+ * \bar{x}_k=b$, the solution set of Equation~\eqref{eq:bipolarFRE} has a greatest element.
\item[(b)] If $a_j^+ * \bar{x}_j\neq b$, for all $j\in\{1,\dots,m\}$, then the set of maximal solutions of Equation~\eqref{eq:bipolarFRE} is finite. Moreover, the number of maximal solutions is:
\[card(\{k\in\{1,\dots,m\}\mid a_k^- * \bar{y}_k=b\})\]
\item[(c)] If there exists some $k\in\{1,\dots,m\}$ such that $a_k^- * \bar{y}_k=b$, the solution set of Equation~\eqref{eq:bipolarFRE} has a least element.
\item[(d)] If $a_j^- *\, \bar{y}_j\neq b$, for all $j\in\{1,\dots,m\}$, then the set of minimal solutions of Equation~\eqref{eq:bipolarFRE} is  finite. Moreover, the number of minimal solutions is:
\[card(\{k\in\{1,\dots,m\}\mid a_k^+ * \bar{x}_k=b\})\]
\end{itemize}}
\end{theorem}
\begin{proof}
(a).
Suppose that there exists $k\in\{1,\dots,m\}$ such that $a_k^+ * \bar{x}_k=b$, and we will demonstrate that the tuple $(\bar{x}_1,\dots,\bar{x}_m)$ is the greatest solution of Equation~\eqref{eq:bipolarFRE}. Clearly, as $(\bar{x}_1,\bar{y}_1,\dots,\bar{x}_m,\bar{y}_m)$ is a solution of Equation~\eqref{eq:FRE}, the inequalities $a_j^+*\bar{x}_j\leq b$ and $a_j^-*\bar{y}_j\leq b$ are verified, for each $j\in\{1,\dots,m\}$. Moreover, because $1\leq \bar{x}_j+ \bar{y}_j$, for all $j\in\{1,\dots,m\}$, and $*$ is an order-preserving operator, we can assert that $a_j^-*(1-\bar{x}_j)\leq b$. Accordingly to the previous inequalities and the fact that $a_k^+ * \bar{x}_k=b$, we conclude that 
\[\bigvee_{j=1}^m (a_j^+* \bar{x}_j) \vee (a_j^-* (1-\bar{x}_j))=b\]
In other words, $(\bar{x}_1,\dots,\bar{x}_m)$ is a solution of Equation~\eqref{eq:bipolarFRE}.

In order to demonstrate that $(\bar{x}_1,\dots,\bar{x}_m)$ is the greatest solution of Equation~\eqref{eq:bipolarFRE}, we will proceed by reductio ad absurdum. Therefore, suppose that there exists a solution $(x_1,\dots,x_m)$ of Equation~\eqref{eq:bipolarFRE} such that $\bar{x}_k<x_k$ for some $k\in\{1,\dots,m\}$, and we will obtain a contradiction. Observe that, since $(x_1,\dots,x_m)$ is a  solution of Equation~\eqref{eq:bipolarFRE}, the tuple $(x_1,1-x_1,\dots,x_m,1-x_m)$ is straightforwardly a solution of Equation~\eqref{eq:FRE}. As a consequence, taking into account that $(\bar{x}_1,\bar{y}_1,\dots,\bar{x}_m,\bar{y}_m)$ is the greatest solution of Equation~\eqref{eq:FRE}, we deduce that, in particular, $x_k\leq \bar{x}_k$. Which is in contradiction with the hypothesis.

Hence, we can conclude that $(\bar{x}_1,\dots,\bar{x}_m)$ is the greatest solution of Equation~\eqref{eq:bipolarFRE}. 

(b). Now, we suppose that $a_j^+ * \bar{x}_j\neq b$ for all $j\in\{1,\dots,m\}$. First of all, notice that we can assert that $a_j^+ * \bar{x}_j< b$, for all $j\in\{1,\dots,m\}$, since $(\bar{x}_1,\bar{y}_1,\dots,\bar{x}_m,\bar{y}_m)$ is a solution of Equation~\eqref{eq:FRE}. Furthermore, as $(\bar{x}_1,\bar{y}_1,\dots,\bar{x}_m,\bar{y}_m)$ is the greatest solution of Equation~\eqref{eq:FRE} and $*$ is continuous in $[0,1]$, we  obtain    that $\bar{x}_j=1$, for all $j\in\{1,\dots,m\}$.

In order to demonstrate Statement (b), consider the set $A=\{k\in\{1,\dots,m\}\mid a_k^- * \bar{y}_k=b\}$, the set $B$ of maximal solutions of Equation~\eqref{eq:bipolarFRE} and the mapping  $f\colon A\to B$ defined as $f(k)=(\bar{x}_1,\dots,\bar{x}_{k-1},1-\bar{y}_k,\bar{x}_{k+1},\dots,\bar{x}_m)$, for all $k\in A$. In the sequel, we will see that $f$ forms a bijection between $A$ and $B$.

Firstly, we will see that $f$ is well-defined. Following an analogous reasoning to the proof in Theorem~\ref{th:luka}, for each $k\in A$, we obtain that the tuple $f(k)=(\bar{x}_1,\dots,\bar{x}_{k-1},1-\bar{y}_k,\bar{x}_{k+1},\dots,\bar{x}_m)$ is a solution of Equation~\eqref{eq:bipolarFRE}. Additionally, it can be easily verified that it is a maximal solution. Specifically, suppose that there exists a solution $(x_1,\dots,x_m)$ of Equation~\eqref{eq:bipolarFRE} such that $f(k)<(x_1,\dots,x_m)$. As $\bar{x}_j=1$ for all $j\in\{1,\dots,m\}$, we have that $x_j=1$, for all $j\in\{1,\dots,m\}$, with $j\neq k$, and $1-\bar{y}_k<x_k$, or equivalently $1-x_k<\bar{y}_k$.

On the one hand, according to the fact that $k\in A$ and due to $*$ is an order-preserving mapping, the following inequality  $a_k^- * (1-x_k)< a_k^- * \bar{y}_k=b$ is satisfied. In addition, since $(x_1,\dots,x_m)$ is a solution of Equation~\eqref{eq:bipolarFRE}, we have that the tuple $(x_1,1-x_1,\dots,x_m,1-x_m)$ is a solution of Equation~\eqref{eq:FRE}. Taking into account that $(\bar{x}_1,\bar{y}_1,\dots,\bar{x}_m,\bar{y}_m)$ is the greatest solution of Equation~\eqref{eq:FRE}, we obtain that $x_k\leq\bar{x}_k$. As a result,
$a_k^+ * x_k\leq a_k^+ * \bar{x}_k<b$ and thus $(a_k^+ * x_k)\vee(a_k^- * (1-x_k))<b$.

On the other hand, for each $j\in\{1,\dots,m\}$ with $j\neq k$, as $x_j=\bar{x}_j=1$ and $a_j^+ * \bar{x}_j<b$ by hypothesis, we can assert that:
\[(a_j^+ * x_j)\vee(a_j^- * (1-x_j))=(a_j^+ * x_j)\vee(a_j^- * 0)=a_j^+ * x_j=a_j^+ * \bar{x}_j<b\]
Hence, since $m$ is finite, we obtain the inequality:
\[\bigvee_{j=1}^m (a_j^+* x_j) \vee (a_j^-* (1-x_j))<b\]
which contradicts the hypothesis of $(x_1,\dots,x_m)$ being a solution of Equation~\eqref{eq:bipolarFRE}.

We can assert then that the mapping $f$ is well-defined. In order to see that $f$ is order-embedding, let $k_1,k_2\in A$ with $k_1\neq k_2$. Notice that, $f(k_1)=f(k_2)$ if and only if $\bar{x}_{k_1}=1-\bar{y}_{k_2}$ and $\bar{x}_{k_2}=1-\bar{y}_{k_1}$. Nevertheless, this is not possible due to $\bar{x}_j=1$ for each $j\in\{1,\dots,m\}$ and $b\neq0$. In particular, the following chain would hold:
\[a_{k_1}^- * \bar{y}_{k_1}=a_{k_1}^- * (1-\bar{x}_{k_2})=a_{k_1}^- * 0=0\neq b\]
Thus, $k_1\notin A$, in contradiction with the hypothesis. Consequently, we conclude that $f(k_1)\neq f(k_2)$ for each $k_1,k_2\in A$ with $k_1\neq k_2$. That is, $f$ is an order-embedding mapping.

To conclude with this demonstration, we will see that $f$ is onto. Given $(\hat{x}_1,\dots,\hat{x}_m)\in B$, that is, a maximal solution of Equation~\eqref{eq:bipolarFRE}, clearly the tuple $(\hat{x}_1,1-\hat{x}_1,\dots,\hat{x}_m,1-\hat{x}_m)$ is a solution of Equation~\eqref{eq:FRE}. Since $(\bar{x}_1,\bar{y}_1,\dots,\bar{x}_m,\bar{y}_m)$ is the greatest solution of Equation~\eqref{eq:FRE}, we obtain that $\hat{x}_j\leq\bar{x}_j$ and  $1-\hat{x}_j\leq\bar{y}_j$ for each $j\in\{1,\dots,m\}$. Therefore, taking into account that $*$ is an order-preserving mapping and $a_j^+*\bar{x}_j<b$ by hypothesis, $a_j^+*\hat{x}_j\leq a_j^+*\bar{x}_j<b$.
As a consequence, we can ensure that there exists $k\in\{1,\dots,m\}$ such that $a_k^-*(1-\hat{x}_k)=b$. Moreover, since $1-\hat{x}_k\leq\bar{y}_k$, the expression $b=a_k^-*(1-\hat{x}_k)\leq a_k^-*\bar{y}_k$ is satisfied. As a result, from the fact that $(\bar{x}_1,\bar{y}_1,\dots,\bar{x}_m,\bar{y}_m)$ is a solution of Equation~\eqref{eq:FRE}, we deduce that $a_k^-*\bar{y}_k=b$, and thus $k\in A$ by definition. In the sequel, we will see that $f(k)=(\hat{x}_1,\dots,\hat{x}_m)$.

Notice that $*$ is a strictly order-preserving mapping and $a_k^-\neq 0$, and therefore the expression $a_k^-*(1-\hat{x}_k)= a_k^-*\bar{y}_k$ implies that $1-\hat{x}_k=\bar{y}_k$. Equivalently, $\hat{x}_k=1-\bar{y}_k$. Hence, in order to prove that $f(k)=(\hat{x}_1,\dots,\hat{x}_m)$, it remains to show that $\hat{x}_j=\bar{x}_j$ for all $j\in\{1,\dots,m\}$, with $j\neq k$. 

Given $j\in\{1,\dots,m\}\setminus\{k\}$, since $a_k^-*(1-\hat{x}_k)=b$ and $\hat{x}_j\leq \bar{x}_j=1$,  the value of $\hat x_j $ in the expression $a_j^-*(1-\hat{x}_j)$ can be as great as possible, maintaining the solvability of the equation.  Thus, we have that $\hat{x}_j=1= \bar{x}_j$.

Hence, we conclude that $f$ forms a bijection between $A$ and $B$, an therefore the number of maximal solution of Equation~\eqref{eq:bipolarFRE} coincides with the cardinal of $A$, as we want to demonstrate.

(c). Assume that there exists some $k\in\{1,\dots,m\}$ such that $a_k^- * \bar{y}_k=b$. Following a dual reasoning to the proof of Statement (a), we can easily demonstrate that the tuple $(1-\bar{y}_1,\dots,1-\bar{y}_m)$ is the least solution of Equation~\eqref{eq:bipolarFRE}.

(d). Now, suppose that $a_j^- *\, \bar{y}_j\neq b$, for each $j\in\{1,\dots,m\}$. Let $A=\{k\in\{1,\dots,m\}\mid a_k^+ * \bar{x}_k=b\}$, let $B$ be the set of maximal solutions of Equation~\eqref{eq:bipolarFRE}, and consider the mapping $f\colon A\to B$ defined as $f(k)=(1-\bar{y}_1,\dots,1-\bar{y}_{k-1},\bar{x}_k,1-\bar{y}_{k+1},\dots,1-\bar{y}_m)$, for each $k\in A$. By a dual reasoning to the proof of Statement (b), we obtain that $f$ forms a bijection between $A$ and $B$, and therefore the number of elements of both sets coincide.
\qed
\end{proof}

{A direct consequence} of the proof of Theorem~\ref{th:luka2} is that we can obtain the analytical description of the maximal and minimal solutions of a solvable bipolar max-product fuzzy equation, as the following corollary shows.

\begin{corollary}\label{cor:luka2_form}
{Let Equation~\eqref{eq:bipolarFRE} be solvable, $b\neq0$ and $(\bar{x}_1,\bar{y}_1,\dots,\bar{x}_m,\bar{y}_m)\in[0,1]^{2m}$ be the greatest solution of its corresponding max-product fuzzy equation~\eqref{eq:FRE}.
\begin{itemize}
\item[(a)] If there exists some $k\in\{1,\dots,m\}$ such that $a_k^+ * \bar{x}_k=b$, then $(\bar{x}_1,\dots,\bar{x}_m)$ is the greatest solution of Equation~\eqref{eq:bipolarFRE}.
\item[(b)] If $a_j^+ * \bar{x}_j\neq b$, for each $j\in\{1,\dots,m\}$, then the set of maximal solutions of Equation~\eqref{eq:bipolarFRE} is finite. The set of maximal solutions of Equation~\eqref{eq:bipolarFRE} is given by:
\[\{(\bar{x}_1,\dots,\bar{x}_{k-1},1-\bar{y}_k,\bar{x}_{k+1},\dots,\bar{x}_m)\mid k\in K_P^-\}\]
where $K_P^-=\{k\in\{1,\dots,m\}\mid a_k^- * \bar{y}_k=b\}$.
\item[(c)] If there exists some $k\in\{1,\dots,m\}$ such that $a_k^- * \bar{y}_k=b$, then $(1-\bar{y}_1,\dots,1-\bar{y}_m)$ is the least solution of Equation~\eqref{eq:bipolarFRE}.
\item[(d)] If $a_j^- *\, \bar{y}_j\neq b$, for each $j\in\{1,\dots,m\}$, then the set of minimal solutions of Equation~\eqref{eq:bipolarFRE} is finite. The set of minimal solutions of Equation~\eqref{eq:bipolarFRE} is given by:
\[\{(1-\bar{y}_1,\dots,1-\bar{y}_{k-1},\bar{x}_k,1-\bar{y}_{k+1},\dots,1-\bar{y}_m)\mid k\in K_P^+\}\]
where $K_P^+=\{k\in\{1,\dots,m\}\mid a_k^+ * \bar{x}_k=b\}$.
\end{itemize}}
\end{corollary}

{The results presented in this section will be illustrated in the following example.}
 
\begin{example}~\label{ex:FE}
We will consider the bipolar max-product fuzzy equation given by Equation~\eqref{ex:bipolarFE}, which contains two unknown variables:
\begin{equation}\label{ex:bipolarFE}
(0.8* x_1) \vee (0.1* (1-x_1))\vee(0.5*x_2)\vee(0.4*(1-x_2))=0.4
\end{equation}
According to the theoretical results presented in~\cite{Medina2016,Medina2017:ija,Perfilieva13}, we obtain that the greatest solution of its corresponding max-product fuzzy equation: 
\begin{equation*}
(0.8* x_1) \vee (0.1* y_1)\vee(0.5*x_2)\vee(0.4*y_2)=0.4
\end{equation*}
is $(0.5,1,0.8,1)$. Taking into account that the inequalities $1\leq 0.5+1$ and $1\leq0.8+1$ are satisfied, Theorem~\ref{th:luka} leads us to assert that Equation~\eqref{ex:bipolarFE} is solvable.

Notice that the equality $0.8*0.5=0.4$ holds. Applying the Statement (a) in Theorem~\ref{th:luka2}, we conclude that Equation~\eqref{ex:bipolarFE} has a greatest solution. Specifically, according to Corollary~\ref{cor:luka2_form}, the greatest solution of Equation~\eqref{ex:bipolarFE} is $(0.5,0.8)$.
Similarly, from the equality $0.4*1=0.4$ and the Statement (c) in Theorem~\ref{th:luka2} and Corollary~\ref{cor:luka2_form}, we obtain that Equation~\eqref{ex:bipolarFE} also has a least solution, which is given by the tuple $(0,0)$. 

Now, we will modify the coefficient which multiplies to the variable $(1-x_2)$ obtaining a new bipolar max-product fuzzy equation:
\begin{equation}\label{ex2:bipolarFRE}
(0.8* x_1) \vee (0.1* (1-x_1))\vee(0.5*x_2)\vee(0.3*(1-x_2))=0.4
\end{equation}
Once again, by using the results included in~\cite{Medina2016,Medina2017:ija,Perfilieva13}, we obtain that the greatest solution of its corresponding max-product fuzzy equation
\begin{equation*}
(0.8* x_1) \vee (0.1* y_1)\vee(0.5*x_2)\vee(0.3*y_2)=0.4
\end{equation*}
is  $(0.5,1,0.8,1)$. Clearly, the hypothesis required in Theorem~\ref{th:luka} are satisfied and hence, Equation~\eqref{ex2:bipolarFRE} is solvable. In this case, we obtain that $0.8*0.5=0.4$, but $0.1*1\neq0.4$ and $0.3*1\neq0.4$. Therefore, taking into account Statements (a) and (d) in Theorem~\ref{th:luka2}, we conclude that Equation~\eqref{ex2:bipolarFRE} has a greatest solution and a finite set of minimal solutions. Specifically, Equation~\eqref{ex2:bipolarFRE} has two minimal solutions since the cardinality of the set $\{j\in\{1,2\}\mid a_j^+ * \bar{x}_j=b\}$ is 2. Corollary~\ref{cor:luka2_form} allows us to conclude that $(0.5,0.8)$ is the greatest solution of {Equation~\eqref{ex2:bipolarFRE}}, whilst  $(0.5,0)$ and  $(0,0.8)$ are its minimal solutions. 
\qed
\end{example}
 
{Observe that, although the tuples $(0.5,0.8)$ and $(0,0)$ are the greatest and the least solution of Equation~\eqref{ex:bipolarFE}, respectively, this does not mean that any tuple between $(0.5,0.8)$ and $(0,0)$ is a solution of Equation~\eqref{ex:bipolarFE}. For instance, $(0.4,0.5)$ is not a solution of Equation~\eqref{ex:bipolarFE}, as shown below
\[(0.8* 0.4) \vee (0.1*(1-0.4))\vee(0.5*0.5)\vee(0.4*(1-0.5))=0.32\neq0.4\]
Therefore, in general,  the solutions of a bipolar max-product fuzzy equation are bounded by the maximal and minimal solutions of the equation, but they are not completely determined by them. The computation of the whole solution set of a bipolar max-product fuzzy equation will be one of the main prospects for future work.}

\section{Bipolar max-product fuzzy relation equations with the standard negation}\label{sec:bipolarFRE}

According to the previous section, one can think that the conditions to ensure the solvability of bipolar max-product fuzzy   equation systems with the standard negation will be the same as the one given in Theorem~\ref{th:luka}. However, this fact is not true in general, as it is shown in the following example. If we consider the next system composed of two max-product fuzzy equations with two unknown variables $x$ and $y$:
\begin{eqnarray*}
(0.5*x)\vee (0.2*y)&=&0.5\\
(0.1*x)\vee (0.8*y)&=&0.4
\end{eqnarray*}
we can easily see that there is only one solution, that is, $(\bar{x},\bar{y})=(1,0.5)$. As a consequence, the unique possible solutions of the system of bipolar max-product fuzzy equations associated with the previous one:
\begin{eqnarray*}
(0.5*x)\vee (0.2*(1-x))&=&0.5\\
(0.1*x)\vee (0.8*(1-x))&=&0.4
\end{eqnarray*}
are $1$ and $(1-0.5)=0.5$. {Nevertheless, the values $1$ and $(1-0.5)=0.5$ are not solutions of this system of bipolar max-product fuzzy equations. Notice that, the value $1$ does not satisfy the second equation and the value $0.5$ does not satisfy the first one, that is:
\begin{eqnarray*}
(0.5*0.5)\vee (0.2*(1-0.5))&=&0.25\vee 0.1=0.25\neq0.5\\
(0.1*1)\vee (0.8*(1-1))&=&0.1\vee 0=0.1\neq0.4
\end{eqnarray*}
Therefore, since the condition $1\leq 1 + 0.5$ is satisfied, we conclude that extra conditions need} to be required in order to guarantee the solvability of a system of bipolar max-product fuzzy equations with the standard negation. These systems are interpreted as a fuzzy relation equation (FRE), as usual.

{From now on, we consider fixed $m,n\in\mathbb{N}$, and $a_{ij}^+,a_{ij}^-,b_i\in[0,1]$, with $i\in \{1,\dots,n\}$ and $j\in\{1,\dots,m\}$.}

\begin{definition}\label{def:FRE}
{Let $x_j\in[0,1]$ be an unknown value, for each $j\in\{1,\dots,m\}$, $\ast$ the product t-norm and $\vee$ the maximum operator. Equation~\eqref{sys:bipolar_general} is called \emph{bipolar max-product fuzzy relation equation with the standard negation}.
\begin{equation}\label{sys:bipolar_general}
      \bigvee_{j=1}^m(a_{ij}^+ \ast x_j) \vee (a_{ij}^-\ast  (1-x_j))=b_i,\qquad i\in\{1,\dots,n\}
  \end{equation}
The \emph{corresponding max-product fuzzy relation equation} of Equation~\eqref{sys:bipolar_general} is given by
\begin{equation}\label{sys:FRE_general}
      \bigvee_{j=1}^m(a_{ij}^+ \ast x_j) \vee (a_{ij}^-\ast  y_j)=b_i,\qquad i\in\{1,\dots,n\}
  \end{equation}
  where $x_j,y_j\in[0,1]$ are unknown values, for each $j\in \{1,\dots,m\}$.
  }
\end{definition}

The next lemma shows two inequalities under which, given a solution of a bipolar max-product FRE, a greater or a smaller solution can be obtained, respectively. {This result will be crucial in the study of the maximal and minimal solutions of an arbitrary bipolar max-product FRE.}

\begin{lemma}\label{lemma:sys_solutionset}
{Let $(x_1,\dots,x_m)$ be a solution of the bipolar max-product FRE~\eqref{sys:bipolar_general} and $(\bar{x}_1,\bar{y}_1,\dots,\bar{x}_m,\bar{y}_m)\in[0,1]^{2m}$ be the greatest solution of its corresponding max-product FRE~\eqref{sys:FRE_general}. The following statements hold:
  \begin{itemize}
  \item[(a)] If there exists $k\in\{1,\dots,m\}$ such that $1-\bar{y}_k<x_k$, then the tuple $(x_1,\dots,x_{k-1},\bar{x}_k,x_{k+1},\dots,x_m)$ is also a solution of Equation~\eqref{sys:bipolar_general}.
  \item[(b)] If there exists $k\in\{1,\dots,m\}$ such that $x_k<\bar{x}_k$, then the tuple $(x_1,\dots,x_{k-1},1-\bar{y}_k,x_{k+1},\dots,x_m)$ is also a solution of Equation~\eqref{sys:bipolar_general}.
  \end{itemize}
  }
\end{lemma}
\begin{proof}
Suppose that there exists $k\in\{1,\dots,m\}$ such that $1-\bar{y}_k<x_k$, and we will prove that Statement (a) holds. Since $(x_1,\dots,x_m)$ is a solution of Equation~\eqref{sys:bipolar_general}, the tuple $(x_1,1-x_1,\dots,x_m,1-x_m)$ is a solution of Equation~\eqref{sys:FRE_general}. Therefore, as $(\bar{x}_1,\bar{y}_1,\dots,\bar{x}_m,\bar{y}_m)$ is the greatest solution of Equation~\eqref{sys:FRE_general}, we can assert that $x_k\leq\bar{x}_k$. Clearly, if $x_k=\bar{x}_k$, then Statement (a) straightforwardly holds. Thus, we assume that $x_k<\bar{x}_k$.  

{Given $i\in\{1,\dots,n\}$, we will prove} that the tuple {$(x_1,\dots,x_{k-1},\bar{x}_k,x_{k+1},\dots,x_m)$} satisfies the $i$-th equation in Equation~\eqref{sys:bipolar_general}. First of all, we consider the trivial case $b_i=0$. 
On the one hand, as $1\leq \bar{x}_k+\bar{y}_k$, and $(\bar{x}_1,\bar{y}_1,\dots,\bar{x}_m,\bar{y}_m)$ is a solution of Equation~\eqref{sys:FRE_general}, we have that
\[
(a_{ik}^+*\bar{x}_k)\vee(a_{ik}^-*(1-\bar x_k))\leq (a_{ik}^+*\bar{x}_k)\vee(a_{ik}^-*\bar{y}_k)=0
\]
On the other hand, since the tuple $(x_1,1-x_1,\dots,x_m,1-x_m)$ is also a solution of Equation~\eqref{sys:FRE_general}, we can assert that:
\[
(a_{ij}^+*x_j)\vee(a_{ij}^-*(1-x_j))=0
\]
for all $j\in\{1,\dots,m\}$. Therefore, $(x_1,\dots,x_{k-1},\bar{x}_k,x_{k+1},\dots,x_m)$ is also a solution of the $i$-th equation in Equation~\eqref{sys:bipolar_general}.

Now, we have that $b_i\neq 0$, that is, $b_i\in \ ]0,1]$. Notice that, if $a_{ik}^+=a_{ik}^-=0$, as $(\bar{x}_1,\bar{y}_1,\dots,\bar{x}_m,\bar{y}_m)$ is a solution of Equation~\eqref{sys:FRE_general}, we obtain that 
\[
(a_{ik}^+*x_k)\vee(a_{ik}^-*(1-x_k))=(a_{ik}^+*\bar{x}_k)\vee(a_{ik}^-*\bar{y}_k)=0< b_i
\]
Otherwise, if $a_{ik}^+\neq 0$ or $a_{ik}^-\neq 0$, taking into account that $*$ is a strictly order-preserving mapping, $(\bar{x}_1,\bar{y}_1,\dots,\bar{x}_m,\bar{y}_m)$ is a solution of Equation~\eqref{sys:FRE_general} and $1-\bar{y}_k<x_k$, or equivalently $1-x_k<\bar{y}_k$, the following chain of inequalities is verified:
\[
(a_{ik}^+*x_k)\vee(a_{ik}^-*(1-x_k))<(a_{ik}^+*\bar{x}_k)\vee(a_{ik}^-*\bar{y}_k)\leq b_i
\]
As a consequence of both cases, since a finite number of variables is considered, $[0,1]$ is totally ordered and $(x_1,\dots,x_m)$ is a solution of Equation~\eqref{sys:bipolar_general}, that is 
\[\bigvee_{j=1}^m(a_{ij}^+*x_j)\vee(a_{ij}^-*(1-x_j))=b_i\]
we can assert that there exists $j'\in\{1,\dots,m\}$ with $j'\neq k$ such that 
\[(a_{ij'}^+*x_{j'})\vee(a_{ij'}^-*(1-x_{j'}))=b_i\]
and clearly, for each $j\in\{1,\dots,m\}$, it is
\[(a_{ij}^+*x_{j})\vee(a_{ij}^-*(1-x_{j}))\leq b_i\]
Hence, in order to prove that the tuple $(x_1,\dots,x_{k-1},\bar{x}_k,x_{k+1},\dots,x_m)$ forms a solution of the $i$-th equation in Equation~\eqref{sys:bipolar_general}, it is sufficient to see that
\[(a_{ij}^+*\bar{x}_k)\vee(a_{ij}^-*(1-\bar{x}_k))\leq b_i\]
But this fact is trivial due to $(\bar{x}_1,\bar{y}_1,\dots,\bar{x}_m,\bar{y}_m)$ is a solution of Equation~\eqref{sys:FRE_general}, and that   $x_k<\bar{x}_k$.  
According to the fact that $i$ is an arbitrary element in $\{1,\dots,n\}$, we conclude that $(x_1,\dots,x_{k-1},\bar{x}_k,x_{k+1},\dots,x_m)$ is a solution of Equation~\eqref{sys:bipolar_general}. Therefore, Statement (a) is proved. Furthermore, following a dual reasoning, Statement (b) is straightforwardly satisfied.
\qed\end{proof}

As a result of this lemma, 
the following corollary can  be easily  obtained.

\begin{corollary}\label{cor:sys_solutions}
{Let $(x_1,\dots,x_m)$ be a solution of the bipolar max-product FRE~\eqref{sys:bipolar_general} and $(\bar{x}_1,\bar{y}_1,\dots,\bar{x}_m,\bar{y}_m)\in[0,1]^{2m}$ be the greatest solution of its corresponding max-product FRE~\eqref{sys:FRE_general}. The following statements hold:
  \begin{itemize}
  \item[(a)] The tuple $(\hat{x}_1,\dots,\hat{x}_m)$ given by $\hat{x}_j=x_j$ if $x_j=1-\bar{y}_j$ and $\hat{x}_j=\bar{x}_j$ otherwise, is also a solution of Equation~\eqref{sys:bipolar_general}.
  \item[(b)] The tuple $(\hat{x}_1,\dots,\hat{x}_m)$ given by $\hat{x}_j=x_j$ if $x_j=\bar{x}_j$ and $\hat{x}_j=1-\bar{y}_j$ otherwise, is also a solution of Equation~\eqref{sys:bipolar_general}.
  \end{itemize}
  }
\end{corollary}
\begin{proof}
Statement (a) (respectively Statement (b)) can {straightforwardly be}  deduced applying Statement (a) (respectively Statement (b)) in Lemma~\ref{lemma:sys_solutionset} for each $k\in\{1,\dots,m\}$ such that $1-\bar{y}_k<x_k$ (respectively $x_k<\bar{x}_k$).
\qed\end{proof}

In order to present a characterization on the resolution of a bipolar max-product FRE, we will introduce the notion of feasible pair of index sets.

\begin{definition}\label{def:feasible}
Let $(\bar{x}_1,\bar{y}_1,\dots,\bar{x}_m,\bar{y}_m)\in[0,1]^{2m}$ be the greatest solution of a solvable max-product FRE~\eqref{sys:FRE_general}. 
  A pair of index sets $(J^+,J^-)\subseteq\{1,\dots,m\}^2$ is said to be \emph{feasible with respect to Equation~\eqref{sys:FRE_general}} if $\bar{x}_j+\bar{y}_j=1$, for each $j\in J^+\cap J^-$, and one of the following statements hold, for each $i\in\{1,\dots,n\}$:
  \begin{enumerate}[label=(\alph*)]
  \item There exists $j\in J^+$ such that $a_{ij}^+*\bar{x}_j=b_i$.
  \item There exists $j\in J^-$ such that $a_{ij}^-*\bar{y}_j=b_i$.
  \end{enumerate}
\end{definition}

The following example illustrates the concept of feasible pair.
{
\begin{example}\label{ex:feasiblepair}
Consider the next max-product FRE 
{\small\begin{equation}\label{eq:FRE_ex}
\begin{array}{ll}
(0.1*x_1) \vee (0.3*y_1) \vee (0.25*x_2) \vee (0.3*y_2) \vee (0.4*x_3) \vee (0.4*y_3)&=0.2\\
(0.1*x_1) \vee (0.8*y_1) \vee (0.3*x_2) \vee (0.6*y_2) \vee (0.3*x_3) \vee (0.5*y_3)&=0.4\\
(0.3*x_1) \vee (0.8*y_1) \vee (0.4*x_2) \vee (0.5*y_2) \vee (0.8*x_3) \vee (0.8*y_3)&=0.4
\end{array}
\end{equation}}
whose greatest solution is the tuple $(\bar{x}_1,\bar{y}_1,\bar{x}_2,\bar{y}_2,\bar{x}_3,\bar{y}_3)=(1,0.5,0.8,0.\widehat{6},0.5,0.5)$. For more details on how to compute the greatest solution of a FRE, see~\cite{Li2008,Medina2017:ija,Peeva2013}.

We will see that the index sets $J^+=\{2\}$ and $J^-=\{1,3\}$ form a feasible pair with respect to Equation~\eqref{eq:FRE_ex}. Observe that, $J^+\cap J^-=\varnothing$.
\begin{itemize}
  \item Case $i=1$ (first equation): we obtain that $2\in J^+$ verifies the equality $a_{12}^+*\bar{x}_2=0.25*0.8=0.2=b_1$.
  \item Case $i=2$ (second equation): we have that $1\in J^-$ satisfies the equality $a_{21}^-*\bar{y}_1=0.8*0.5=0.4=b_2$.
  \item Case $i=3$ (third equation): we obtain that $3\in J^-$ verifies the equality $a_{33}^-*\bar{y}_3=0.8*0.5=0.4=b_3$.
\end{itemize}
Hence, the pair $(J^+,J^-)$ forms a feasible pair with respect to Equation~\eqref{eq:FRE_ex}.
\qed
\end{example}
}

{The reader may realize that the intuition behind the notion of feasible pair is unquestionably close to the concept of covering~\cite{Lin2011,Markovskii}. Indeed, the concept of feasible pair follows the philosophy of the covering problem for the case of bipolar fuzzy relation equations. Therefore, the definition of feasible pair can be seen as a modification of the concept of covering in order to consider the bipolar framework.}

The proposition below presents {a useful property} of feasible pairs. Specifically, given a feasible pair $(J^+,J^-)$, we can obtain new feasible pairs adding to, either $J^+$ or $J^-$, any index $k\in\{1,\dots, m\}$ satisfying the equality $\bar{x}_k+\bar{y}_k=1$.

\begin{proposition}\label{prop:feasible}
Let $m\in\mathbb{N}$, $(\bar{x}_1,\bar{y}_1,\dots,\bar{x}_m,\bar{y}_m)\in[0,1]^{2m}$ be the greatest solution of a solvable max-product FRE~\eqref{sys:FRE_general}, $(J^+,J^-)$ be a feasible pair with respect to Equation~\eqref{sys:FRE_general} and $k\in\{1,\dots,m\}$ such that $\bar{x}_k+\bar{y}_k=1$. Then $(J^+\cup\{k\},J^-)$ and $(J^+,J^-\cup\{k\})$ are feasible pairs with respect to Equation~\eqref{sys:FRE_general}.
\end{proposition}
\begin{proof}
We will show that $(J^+\cup\{k\},J^-)$ is a feasible pair with respect to Equation~\eqref{sys:FRE_general}. The proof for $(J^+,J^-\cup\{k\})$ is analogous.

On the one hand, we consider $j\in (J^+\cup\{k\})\cap J^-$. If $j\neq k$, we obtain that $j\in J^+\cap J^-$. Therefore, since $(J^+,J^-)$ is a feasible pair, the equality $\bar{x}_j+\bar{y}_j=1$ holds. Furthermore, if $j=k$, then by hypothesis $\bar{x}_k+\bar{y}_k=1$.

On the other hand, taking into account that $(J^+,J^-)$ is a feasible pair with respect to Equation~\eqref{sys:FRE_general} and $J^+\subseteq J^+\cup\{k\}$, we obtain that, for each $i\in\{1,\dots,n\}$, one the following statements is satisfied:
\begin{enumerate}[label=(\alph*)]
  \item There exists $j\in J^+\cup\{k\}$ such that $a_{ij}^+*\bar{x}_j=b_i$.
  \item There exists $j\in J^-$ such that $a_{ij}^-*\bar{y}_j=b_i$.
  \end{enumerate}
Hence, we conclude that $(J^+\cup\{k\},J^-)$ is a feasible pair with respect to Equation~\eqref{sys:FRE_general}.
\qed\end{proof}

In what follows, we introduce a result which shows under what conditions an arbitrary bipolar max-product fuzzy relation equation with the standard negation is solvable. In particular, the solvability of bipolar max-product FREs is characterized by the existence of feasible pairs. {The idea behind this characterization is outlined in the sequel.

Notice that, in Example~\ref{ex:feasiblepair}, the pair $(J^+,J^-)$, where $J^+=\{2\}$ and $J^-=\{1,3\}$, is a feasible pair with respect to Equation~\eqref{eq:FRE_ex}. The underlying significance of this fact is that all the independent terms of Equation~\eqref{eq:FRE_ex} can be reached by means of $\bar{x}_2$ (since $J^+=\{2\}$), $\bar{y}_1$ and $\bar{y}_3$ (since $J^-=\{1,3\}$). Hence, the existence of a feasible pair is in a certain way a sufficient condition for the solvability of a bipolar max-product FRE, whenever each subequation of the bipolar equation is solvable. In other words, as long as $1\leq\bar{x}_1+\bar{y}_1$ and $1\leq\bar{x}_2+\bar{y}_2$. As we will formalize below, the existence of a feasible pair together with the inequality $1\leq\bar{x}_j+\bar{y}_j$, for each $j\in\{1,\dots,m\}$, is not only a sufficient condition, but a necessary condition for the solvability of a bipolar max-product FRE.}

\begin{theorem}\label{th:sys_general_mxn}
{The bipolar max-product  FRE~\eqref{sys:bipolar_general} is solvable if and only if there exists (at least) a feasible pair $(J^+,J^-)$ with respect to its corresponding max-product FRE~\eqref{sys:FRE_general} and the inequality $1\leq\bar{x}_j+\bar{y}_j$ holds for each $j\in\{1,\dots,m\}$, where $(\bar{x}_1,\bar{y}_1,\dots,\bar{x}_m,\bar{y}_m)\in[0,1]^{2m}$ is the greatest solution of Equation~\eqref{sys:FRE_general}.
}
\end{theorem}
\begin{proof}
To begin with, suppose that $1\leq\bar{x}_j+\bar{y}_j$ holds, for each $j\in\{1,\dots,m\}$, and  that there exists  a feasible pair $(J^+,J^-)$  with respect to Equation~\eqref{sys:FRE_general}.
Consider the tuple $(\hat{x}_1,\dots,\hat{x}_m)$ defined, for each $j\in\{1,\dots,m\}$, as:
\[\hat{x}_j= \left\{\begin{array}{ll}
             \bar{x}_j & \hbox{ if}\quad j\in J^+ \\
             1-\bar{y}_j & \hbox{ otherwise}
             \end{array}
  \right.\]
We will see that $(\hat{x}_1,\dots,\hat{x}_m)$ is a solution of Equation~\eqref{sys:bipolar_general}. Consider $i\in\{1,\dots,n\}$ fixed, by hypothesis, one of the following holds:
  \begin{enumerate}[label=(\alph*)]
  \item There exists $j\in J^+$ such that $a_{ij}^+*\bar{x}_j=b_i$.
  \item There exists $j\in J^-$ such that $a_{ij}^-*\bar{y}_j=b_i$.
  \end{enumerate}
If Statement (a) is verified, since $j\in J^+$, we obtain that $\hat{x}_j=\bar{x}_j$, and thus $a_{ij}^+*\hat{x}_j=b_i$. In addition, according to the fact that $1\leq\bar{x}_j+\bar{y}_j$, the operator $*$ is order-preserving and the tuple $(\bar{x}_1,\bar{y}_1,\dots,\bar{x}_m,\bar{y}_m)\in[0,1]^{2m}$ is a solution of Equation~\eqref{sys:FRE_general}, we can assert that the following chain of inequalities holds:
\[a_{ij}^-*(1-\hat{x}_j)\leq a_{ij}^-*\bar{y}_j\leq b_i\]
On the contrary, if Statement (b) is satisfied, there exists $j\in J^-$ such that $a_{ij}^-*\bar{y}_j=b_i$. As a result, $\hat{x}_j$ is defined as $1-\bar{y}_j$, and therefore:
\[a_{ij}^-*(1-\hat{x}_j)=a_{ij}^-*\bar{y}_j=b_i\]
Furthermore, as $1\leq\bar{x}_j+\bar{y}_j$, $*$ is order-preserving and $(\bar{x}_1,\bar{y}_1,\dots,\bar{x}_m,\bar{y}_m)\in[0,1]^{2m}$ is a solution of Equation~\eqref{sys:FRE_general}, we obtain that:
\[a_{ij}^+*\hat{x}_j= a_{ij}^+*(1-\bar{y}_j)\leq a_{ij}^+*\bar{x}_j\leq b_i\]
Hence, in both cases, we can ensure that there exists $j\in\{1,\dots,m\}$ such that: 
\[(a_{ij}^+ \ast \hat{x}_j) \vee (a_{ij}^-\ast  (1-\hat{x}_j))=b_i\]
Now, we will see that, for any $j^*\in\{1,\dots,m\}$ different from $j$, we obtain that:
\[(a_{i{j^*}}^+ \ast \hat{x}_{j^*}) \vee (a_{i{j^*}}^-\ast  (1-\hat{x}_{j^*}))\leq b_i\]
On the one hand, if $j^*\in J^+$, $\hat{x}_{j^*}$ coincides with $\bar{x}_{j^*}$. Therefore, taking into account that $1\leq\bar{x}_{j^*}+\bar{y}_{j^*}$ and that $(\bar{x}_1,\bar{y}_1,\dots,\bar{x}_m,\bar{y}_m)\in[0,1]^{2m}$ is a solution of Equation~\eqref{sys:FRE_general}, we can assert that: 
\[(a_{i{j^*}}^+ \ast \hat{x}_{j^*}) \vee (a_{i{j^*}}^-\ast  (1-\hat{x}_{j^*}))\leq b_i\]
On the other hand, if $j^*\notin J^+$, then $\hat{x}_{j^*}=1-\bar{y}_{j^*}$. Due to $1\leq\bar{x}_{j^*}+\bar{y}_{j^*}$, $*$ is order-preserving and $(\bar{x}_1,\bar{y}_1,\dots,\bar{x}_m,\bar{y}_m)\in[0,1]^{2m}$ is a solution of Equation~\eqref{sys:FRE_general}, we deduce the next chain of inequalities
\[(a_{i{j^*}}^+ \ast \hat{x}_{j^*}) \vee (a_{i{j^*}}^-\ast  (1-\hat{x}_{j^*}))=(a_{i{j^*}}^+ \ast (1-\bar{y}_{j^*})) \vee (a_{i{j^*}}^-\ast  \bar{y}_{j^*})\leq (a_{i{j^*}}^+*\bar{x}_{j^*})\vee (a_{i{j^*}}^-*\bar{y}_{j^*})\leq b_i\]
This leads us to conclude that
\[\bigvee_{j=1}^m(a_{ij}^+ \ast \hat{x}_j) \vee (a_{ij}^-\ast  (1-\hat{x}_j))=b_i\]
that is, $(\hat{x}_1,\dots,\hat{x}_m)$ is a solution of equation $i$ in Equation~\eqref{sys:bipolar_general}. According to the fact that $i$ is an arbitrary element in $\{1,\dots,n\}$, we conclude that $(\hat{x}_1,\dots,\hat{x}_m)$ is a solution of Equation~\eqref{sys:bipolar_general}, {as we want to prove}.

Now, suppose that Equation~\eqref{sys:bipolar_general} is solvable and that $(\hat{x}_1,\dots,\hat{x}_m)$ is a solution of Equation~\eqref{sys:bipolar_general}. Clearly, the tuple $(\hat{x}_1,1-\hat{x}_1,\dots,\hat{x}_m,1-\hat{x}_m)$ is a solution of Equation~\eqref{sys:FRE_general}. Furthermore, according to the fact that $(\bar{x}_1,\bar{y}_1,\dots,\bar{x}_m,\bar{y}_m)$ is the greatest solution of Equation~\eqref{sys:FRE_general}, the inequalities $\hat{x}_j\leq\bar{x}_j$ and $1-\hat{x}_j\leq\bar{y}_j$ hold, for each $j\in\{1,\dots,m\}$. Therefore, we obtain that $1\leq\hat{x}_j+\bar{y}_j\leq \bar{x}_j+\bar{y}_j$, for each $j\in\{1,\dots,m\}$.

According to Corollary~\ref{cor:sys_solutions}, we can suppose without loss of generality that, for each $j\in\{1,\dots,m\}$, either $\hat{x}_j=\bar{x}_j$ or $\hat{x}_j=1-\bar{y}_j$. Consider the two index sets $J^+,J^-\subseteq\{1,\dots,m\}$ defined as follows:
\begin{eqnarray*}
J^+&=&\{j\in\{1,\dots,m\}\mid \hat{x}_j=\bar{x}_j\}\\
J^-&=&\{j\in\{1,\dots,m\}\mid \hat{x}_j=1-\bar{y}_j\}
\end{eqnarray*}

In the sequel, we will prove that $(J^+,J^-)$ forms a feasible pair with respect to Equation~\eqref{sys:FRE_general}. That is, we will see that $1=\bar{x}_j+\bar{y}_j$ for each $j\in J^+\cap J^-$, and either Statement (a) or (b) is satisfied for each $i\in\{1,\dots,n\}$. In fact, the former is straightforwardly obtained since, if $j\in J^+\cap J^-$, then $\bar{x}_j=1-\bar{y}_j$. Consequently, $1=\bar{x}_j+\bar{y}_j$.

{Given $i\in\{1,\dots,n\}$}, as $(\hat{x}_1,\dots,\hat{x}_m)$ is a solution of Equation~\eqref{sys:bipolar_general}, we can assert that one of the following statements hold:
\begin{itemize}
\item there exists $j\in\{1,\dots,m\}$ such that $a_{ij}^+*\hat{x}_j=b_i$
\item there exists $j\in\{1,\dots,m\}$ such that $a_{ij}^-*(1-\hat{x}_j)=b_i$
\end{itemize}

Assume that there exists $j\in\{1,\dots,m\}$ such that $a_{ij}^+*\hat{x}_j=b_i$. Clearly, if $\hat{x}_j=\bar{x}_j$, then $j\in J^+$ and Statement (a) of Definition~\ref{def:feasible} is straightforwardly satisfied. 

 {Now, we consider  that $\hat{x}_j\neq\bar{x}_j$, then, by Corollary~\ref{cor:sys_solutions} and the supposition above, we have that $\hat{x}_j=1-\bar{y}_j$,
 and so, $j\in J^-$.} {Notice} that, since $(\bar{x}_1,\bar{y}_1,\dots,\bar{x}_m,\bar{y}_m)$ is a solution of Equation~\eqref{sys:FRE_general}, then $a_{ij}^+*\bar{x}_j\leq b_i$ and $a_{ij}^-*\bar{y}_j\leq b_i$, for each $j\in\{1,\dots,m\}$. Furthermore, taking into account that $*$ is an order-preserving mapping and $(\hat{x}_1,1-\hat{x}_1,\dots,\hat{x}_m,1-\hat{x}_m)\leq(\bar{x}_1,\bar{y}_1,\dots,\bar{x}_m,\bar{y}_m)$, we obtain that
$a_{ij}^-*\hat{x}_j\leq a_{ij}^-*\bar{x}_j\leq b_i$, for each $j\in\{1,\dots,m\}$. Therefore, since $a_{ij}^+*\hat{x}_j=b_i$ holds, the equality $a_{ij}^+*\bar{x}_j=b_i$ is also satisfied. As $\hat{x}_j\neq\bar{x}_j$ and $*$ is strictly order-preserving, {we deduce   that $a_{ij}^+=0$ and so, $b_i=0$}. Hence, as $a_{ij}^-*(1-\hat{x}_j)\leq 0$, 
we conclude that $a_{ij}^-*(1-\hat{x}_j)=0$. In other words, Statement (b) holds.

On the contrary, assume that there exists $j\in\{1,\dots,m\}$ such that $a_{ij}^-*(1-\hat{x}_j)=b_i$. Clearly, if $\hat{x}_j=1-\bar{y}_j$, then $j\in J^-$ and Statement (b) of Definition~\ref{def:feasible} is straightforwardly satisfied. 
Following an analogous reasoning to the previous case, we deduce that, if {$1-\hat{x}_j\neq\bar{y}_j$, then  $j\in J^+$ and  Statement (a) of Definition~\ref{def:feasible} holds}. As a result, we obtain that $(J^+,J^-)$ forms a feasible pair.
\qed\end{proof}

Hence, this theorem presents necessary and sufficient conditions in order to guarantee the solvability of a bipolar max-product FRE. {Moreover, the proof of Theorem~\ref{th:sys_general_mxn} shows how to define different solutions of this equation. The next section will be focused on deepening in the study bipolar max-product FREs, characterizing the maximal and minimal solutions.

\section{{Characterizing the maximal and minimal solutions of bipolar max-product FREs}}\label{sec:allsol}

This section will study whether the existence of maximal or minimal solutions is always guaranteed for bipolar max-product FREs and, if they exist, how they can be computed. In what follows, {we will achieve these goals} introducing the conditions under which an arbitrary bipolar max-product FRE has a greatest solution (respectively a least solution) or maximal solutions (respectively minimal solutions).

The idea behind the existence theorem of the greatest solution (respectively a least solution) or maximal solutions (respectively minimal solutions) of a bipolar max-product FRE is quite different from the case of one bipolar max-product fuzzy equation. 

Hence, we will consider the set $S^+$ of index sets that appear in the first argument of feasible pairs in order to  provide a characterization of its greatest solution or maximal solutions. Specifically, we obtain that the maximal solutions of such set $S^+${, with respect to the inclusion $\subseteq$, are closely related to the maximal solutions of the bipolar fuzzy relation equation.}

\begin{theorem}\label{th:sys_maximal}
{Let Equation~\eqref{sys:bipolar_general} be a solvable bipolar max-product FRE, $(\bar{x}_1,\bar{y}_1,\dots,\bar{x}_m,\bar{y}_m)\in[0,1]^{2m}$ be the greatest solution of its corresponding max-product FRE~\eqref{sys:FRE_general}, $S$ be the set of feasible pairs with respect to Equation~\eqref{sys:FRE_general} and $S^+=\{J^+\mid(J^+,J^-)\in S\}$. The following statements hold:
  \begin{enumerate}
  \item[(a)] If the set $S^+$ has a greatest element, then Equation~\eqref{sys:bipolar_general} has a greatest solution.
  \item[(b)] The number of maximal solutions of Equation~\eqref{sys:bipolar_general} coincides with the number of maximal elements of $S^+$.
  \end{enumerate}
  }
\end{theorem}
\begin{proof}
In order to demonstrate Statement (a), suppose that there exists the greatest element $\hat{J}^+$ of $S^+$ and consider the tuple $(\hat{x}_1,\dots,\hat{x}_m)$ given by
  \[\hat{x}_j= \left\{\begin{array}{ll}
             \bar{x}_j & \hbox{ if}\quad j\in \hat{J}^+ \\
             1-\bar{y}_j & \hbox{ otherwise}
             \end{array}
  \right.\]
We will see that $(\hat{x}_1,\dots,\hat{x}_m)$ is the greatest solution of Equation~\eqref{sys:bipolar_general}. Notice that, following a similar reasoning to the proof of Theorem~\ref{th:sys_general_mxn}, we can easily obtain that $(\hat{x}_1,\dots,\hat{x}_m)$ is a solution of Equation~\eqref{sys:bipolar_general}. It remains to see then that it is its greatest solution. We will proceed by reductio ad absurdum.

Suppose that there exists a solution $(x_1,\dots,x_m)$ of Equation~\eqref{sys:bipolar_general} such that $(x_1,\dots,x_m)\not\leq(\hat{x}_1,\dots,\hat{x}_m)$. As a result, there exists $k\in\{1,\dots,m\}$ such that $\hat{x}_k<x_k$. Notice that, the tuple $(x_1,1-x_1,\dots,x_m,1-x_m)$ forms a solution of Equation~\eqref{sys:FRE_general}, and thus $x_j\leq\bar{x}_j$ for each $j\in\{1,\dots,m\}$. In particular, $\hat{x}_k<x_k\leq\bar{x}_k$. Therefore, according to the definition of $(\hat{x}_1,\dots,\hat{x}_m)$, we can assert that $k\notin\hat{J}^+$, and thus $1-\bar{y}_k=\hat{x}_k<x_k$. Furthermore, taking into account Statement (a) in Lemma~\ref{lemma:sys_solutionset}, we obtain that the tuple $(x_1,\dots,x_{k-1},\bar{x}_k,x_{k+1},\dots,x_m)$ is a solution of Equation~\eqref{sys:bipolar_general}. In order to avoid including a new notation, we will suppose without loss of generality that $(x_1,\dots,x_m)=(x_1,\dots,x_{k-1},\bar{x}_k,x_{k+1},\dots,x_m)$, that is, $x_k=\bar{x}_k$.

Hence, consider the index sets $J^+$ and $J^-$ defined as follows
\begin{eqnarray*}
J^+&=&\{j\in\{1,\dots,m\}\mid x_j=\bar{x}_j\}\\
J^-&=&\{j\in\{1,\dots,m\}\mid x_j=1-\bar{y}_j\}
\end{eqnarray*}
Clearly, following an analogous reasoning to the proof in Theorem~\ref{th:sys_general_mxn}, we obtain that $(J^+,J^-)\in S$. Besides, $k\in J^+$, and thus we can assert that $J^+\not\subseteq\hat{J}^+$, which contradicts the hypothesis of $\hat{J}^+$ being the greatest element of $S^+$.

As a result, we conclude that $(\hat{x}_1,\dots,\hat{x}_m)$ is the greatest solution of Equation~\eqref{sys:bipolar_general}, as we want to demonstrate, and thus Statement (a) holds.

In what regards Statement (b), let $A$ be the set of maximal elements of $S^+$, $B$ be the set of maximal solutions of Equation~\eqref{sys:bipolar_general} and consider the mapping $f\colon A\to B$ which maps each maximal element $\hat{J}^+$ of $S^+$ with the tuple $(\hat{x}_1,\dots,\hat{x}_m)$ given by
\[\hat{x}_j= \left\{\begin{array}{ll}
             \bar{x}_j & \hbox{ if}\quad j\in \hat{J}^+ \\
             1-\bar{y}_j & \hbox{ otherwise}
             \end{array}
  \right.\]
In the following, we will prove that $f$ is a bijection. Firstly, let us see that $f$ is well-defined. Clearly, by the proof in Theorem~\ref{th:sys_general_mxn}, given $\hat{J}^+\in S^+$, the tuple $f(\hat{J}^+)=(\hat{x}_1,\dots,\hat{x}_m)$ is a solution of Equation~\eqref{sys:bipolar_general}. Therefore, it remains to see that it is a maximal solution of Equation~\eqref{sys:bipolar_general}. We will proceed by reductio ad absurdum.

Suppose that $(\hat{x}_1,\dots,\hat{x}_m)$ is not a maximal solution of Equation~\eqref{sys:bipolar_general}, that is, there exists a solution $(x_1,\dots,x_m)$ of Equation~\eqref{sys:bipolar_general} such that $(\hat{x}_1,\dots,\hat{x}_m)<(x_1,\dots,x_m)$. Consider the sets $J^+$ and $J^-$ given by
\begin{eqnarray*}
J^+&=&\{j\in\{1,\dots,m\}\mid x_j=\bar{x}_j\}\\
J^-&=&\{j\in\{1,\dots,m\}\mid x_j=1-\bar{y}_j\}
\end{eqnarray*}
By an analogous reasoning to the proof in Theorem~\ref{th:sys_general_mxn}, $(J^+,J^-)\in S$. Now, observe that $(x_1,1-x_1,\dots,x_m,1-x_m)$ is a solution of Equation~\eqref{sys:FRE_general}, and thus $x_j\leq\bar{x}_j$ for each $j\in\{1,\dots,m\}$. This fact together with the inequality $(\hat{x}_1,\dots,\hat{x}_m)<(x_1,\dots,x_m)$ allows us to assert that $x_j=\bar{x}_j$, for each $j\in\hat{J}^+$. As a result, $\hat{J}^+\subseteq J^+$. Furthermore, as $(\hat{x}_1,\dots,\hat{x}_m)$ is strictly smaller than $(x_1,\dots,x_m)$, we deduce that there exists $k\in\{1,\dots,m\}$ such that $\hat{x}_k<x_k$. Clearly, it must be $\hat{x}_k=1-\bar{y}_k$ and thus, by definition, $k\notin\hat{J}^+$. Lastly, we have just to realise that, basing on Statement (a) in Lemma~\ref{lemma:sys_solutionset}, we can suppose without loss of generality that $x_k=\bar{x}_k$. This fact implies that $k\in J^+$, and thus $\hat{J}^+\subset J^+$, in contradiction with the hypothesis. Hence, we conclude that the mapping $f$ is well-defined.

In order to see that $f$ is order-embedding, given $\hat{J}^+_1,\hat{J}^+_2\in A$ with $\hat{J}^+_1\neq\hat{J}^+_2$, we suppose that $f(\hat{J}^+_1)=f(\hat{J}^+_2)$ and we will obtain a contradiction. In fact, if $\hat{J}^+_1\neq\hat{J}^+_2$, then we can suppose without loss of generality that there exists $k\in\hat{J}^+_1$ such that $k\notin\hat{J}^+_2$. As a consequence, given $f(\hat{J}^+_1)=(\hat{x}_1^1,\dots,\hat{x}_m^1)$ and $f(\hat{J}^+_2)=(\hat{x}_1^2,\dots,\hat{x}_m^2)$, then by definition $\hat{x}_k^1=\bar{x}_k$ and $\hat{x}_k^2=1-\bar{y}_k$. Therefore, if the equality $f(\hat{J}^+_1)=f(\hat{J}^+_2)$ is satisfied, we obtain that $\bar{x}_k=1-\bar{y}_k$. In other words, $1=\bar{x}_k+\bar{y}_k$.

Now, as $\hat{J}^+_2\in S^+$, notice that there exists some index set $\hat{J}^-_2\subseteq\{1,\dots,m\}$ such that $(\hat{J}^+_2,\hat{J}^-_2)\in S$. Furthermore, consider the sets $J^+=\hat{J}^+_2\cup\{k\}$ and $J^-=\hat{J}^-_2$. Notice that, since $1=\bar{x}_k+\bar{y}_k$, Proposition~\ref{prop:feasible} allows us to ensure that $(J^+,J^-)$ is also a feasible pair with respect to Equation~\eqref{sys:FRE_general}, that is, $(J^+,J^-)\in S$. Hence, $J^+\in S^+$ and clearly $\hat{J}^+_2\subset J^+$, in contradiction with the hypothesis of $\hat{J}^+_2$ being a maximal element of $S^+$. We can conclude then that $f$ is order-embedding.

It remains to demonstrate that $f$ is onto. Given $(\hat{x}_1,\dots,\hat{x}_m)\in B$, consider the sets $\hat{J}^+$ and $\hat{J}^-$ defined as follows
\begin{eqnarray*}
\hat{J}^+&=&\{j\in\{1,\dots,m\}\mid \hat{x}_j=\bar{x}_j\}\\
\hat{J}^-&=&\{j\in\{1,\dots,m\}\mid \hat{x}_j=1-\bar{y}_j\}
\end{eqnarray*}
In the following, we will see that $\hat{J}^+\in A$ and $f(\hat{J}^+)=(\hat{x}_1,\dots,\hat{x}_m)$. It is important to highlight that, if there exists $j\in\{1,\dots,m\}$ such that $1-\bar{y}_j<\hat{x}_j<\bar{x}_j$, then Lemma~\ref{lemma:sys_solutionset} leads us to assert that there exists a solution of Equation~\eqref{sys:bipolar_general} that is greater than $(\hat{x}_1,\dots,\hat{x}_m)$, and therefore $(\hat{x}_1,\dots,\hat{x}_m)\notin B$, in contradiction with the hypothesis. Moreover, the cases $\bar{x}_j<\hat{x}_j$ and $\hat{x}_j <1-\bar{y}_j$ are not possible since $(\bar{x}_1,\bar{y}_1,\dots,\bar{x}_m,\bar{y}_m)$ is the greatest solution of Equation~\eqref{sys:FRE_general}. As a consequence, we can suppose that, for each $j\in\{1,\dots,m\}$, either $\hat{x}_j=\bar{x}_j$ or $\hat{x}_j=1-\bar{y}_j$. On the one hand, this fact implies that, whenever $\hat{J}^+\in A$, we can ensure that $f(\hat{J}^+)=(\hat{x}_1,\dots,\hat{x}_m)$. On the other hand, we obtain that the expression $\hat{J}^+\cup\hat{J}^-=\{1,\dots,m\}$ holds. In the following, to finish with this demonstration, we will see that indeed $\hat{J}^+\in A$.

By an analogous reasoning to the proof of Theorem~\ref{th:sys_general_mxn}, we obtain that $(\hat{J}^+,\hat{J}^-)\in S$, and thus $\hat{J}^+\in S^+$. Suppose now that $\hat{J}^+\notin A$, that is, that there exists $J^+\in S^+$ such that $\hat{J}^+\subset J^+$, and consider the tuple $(x_1,\dots,x_m)$ defined as
\[x_j= \left\{\begin{array}{ll}
             \bar{x}_j & \hbox{ if}\quad j\in J^+ \\
             1-\bar{y}_j & \hbox{ otherwise}
             \end{array}
  \right.\]
Clearly, following an analogous reasoning to the proof of Theorem~\ref{th:sys_general_mxn}, $(x_1,\dots,x_m)$ is a solution of Equation~\eqref{sys:bipolar_general}. By definition, for each $j\in\hat{J}^+\subset J^+$, we obtain that $\hat{x}_j=\bar{x}_j=x_j$. In addition, there exists $k\notin\hat{J}^+$ such that $k\in J^+$, which implies that $\hat{x}_k\neq \bar{x}_k$ and $x_k=\bar{x}_k$. Therefore, taking into account that $(\hat{x}_1,1-\hat{x}_1,\dots,\hat{x}_m,1-\hat{x}_m)$ is a solution of Equation~\eqref{sys:FRE_general} and consequently $\hat{x}_k\leq\bar{x}_k$, we can assert that $\hat{x}_k< \bar{x}_k=x_k$ for this index $k$. Finally, for each $j\notin \hat{J}^+$, as $\hat{J}^+\cup\hat{J}^-=\{1,\dots,m\}$, we obtain that $j\in \hat{J}^-$, from which $\hat{x}_j=1-\bar{y}_j$. Taking into account that $(\hat{J}^+,\hat{J}^-)\in S$, we can assert that $1\leq \bar{x}_j+\bar{y}_j$, or equivalently $1-\bar{y}_j\leq \bar{x}_j$. As a result, the inequality $\hat{x}_j\leq x_j$ holds.

Consequently, we deduce that $(\hat{x}_1,\dots,\hat{x}_m)<(x_1,\dots,x_m)$, which leads us to a contradiction due to $(\hat{x}_1,\dots,\hat{x}_m)\in B$ by hypothesis. Hence, we can assert that $\hat{J}^+\in A$.
\qed\end{proof}

Theorem~\ref{th:sys_maximal} does not only provides a characterization of the existence of the greatest solution or maximal solutions of a solvable bipolar max-product fuzzy relation equation. From its demonstration, we can also deduce the analytical description of these solutions which are provided in the next corollary.

\begin{corollary}\label{cor:sys_maximal}
{Let Equation~\eqref{sys:bipolar_general} be a solvable bipolar max-product FRE, $(\bar{x}_1,\bar{y}_1,\dots,\bar{x}_m,\bar{y}_m)\in[0,1]^{2m}$ be the greatest solution of its corresponding max-product FRE~\eqref{sys:FRE_general}, $S$  be the set of feasible pairs with respect to Equation~\eqref{sys:FRE_general} and $S^+=\{J^+\mid(J^+,J^-)\in S\}$. Consider the mapping $f\colon S^+\to [0,1]^m$ which associates each $J^+\in S^+$ with the tuple $(x_1,\dots,x_m)$ defined, for each $j\in\{1,\dots,m\}$, as:
\[x_j= \left\{\begin{array}{ll}
             \bar{x}_j & \hbox{ if}\quad j\in J^+ \\
             1-\bar{y}_j & \hbox{ otherwise}
             \end{array}
  \right.\]
The following statements hold:
  \begin{enumerate}
  \item[(a)] If the set $S^+$ has a greatest element $J^+$, then $f(J^+)$ is the greatest solution of Equation~\eqref{sys:bipolar_general}.
  \item[(b)] Let $M^+$ be the set of maximal elements of $S^+$. Then, the set of maximal solutions of Equation~\eqref{sys:bipolar_general} is given by:
  \[\{f(J^+)\mid J^+\in M^+\}\]
  \end{enumerate}
  }
\end{corollary}

In what regards the existence of the least solution or minimal solutions of a given bipolar max-product FRE, the idea is dual to the previous case. Now, we will consider the set composed of index sets appearing in the second argument of feasible pairs with respect to its corresponding max-product fuzzy equation.
\begin{theorem}\label{th:sys_minimal}
{Let Equation~\eqref{sys:bipolar_general} be a solvable bipolar max-product FRE, $(\bar{x}_1,\bar{y}_1,\dots,\bar{x}_m,\bar{y}_m)\in[0,1]^{2m}$ be the greatest solution of its corresponding max-product FRE~\eqref{sys:FRE_general}, $S$  be the set of feasible pairs with respect to Equation~\eqref{sys:FRE_general} and $S^-=\{J^-\mid(J^+,J^-)\in S\}$. The following statements hold:
  \begin{enumerate}
  \item[(a)] If the set $S^-$ has a greatest element, then Equation~\eqref{sys:bipolar_general} has a least solution.
  \item[(b)] The number of minimal solutions of Equation~\eqref{sys:bipolar_general} coincides with the number of maximal elements of $S^-$.
  \end{enumerate}
  }
\end{theorem}
\begin{proof}
As far as Statement (a) is concerned, suppose that there exists the greatest element $\hat{J}^-$ of $S^-$ and consider the tuple $(\hat{x}_1,\dots,\hat{x}_m)$ given by
  \[\hat{x}_j= \left\{\begin{array}{ll}
             1-\bar{y}_j & \hbox{ if}\quad j\in \hat{J}^- \\
             \bar{x}_j & \hbox{ otherwise}
             \end{array}
  \right.\]
Following a dual reasoning to the proof of Theorem~\ref{th:sys_maximal} and taking into account Statement (b) in Lemma~\ref{lemma:sys_solutionset}, we deduce that $(\hat{x}_1,\dots,\hat{x}_m)$ is the least solution of Equation~\eqref{sys:bipolar_general}. 

Regarding Statement (b), we consider the set  $A$ of maximal elements of $S^-$,   the set  $B$  of minimal solutions of Equation~\eqref{sys:bipolar_general} and   the mapping $f\colon A\to B$, which maps each maximal element $\hat{J}^-$ of $S^-$ to the tuple $(\hat{x}_1,\dots,\hat{x}_m)$, defined for each $j\in\{1,\dots,m\}$ as
\[\hat{x}_j= \left\{\begin{array}{ll}
             1-\bar{y}_j & \hbox{ if}\quad j\in \hat{J}^- \\
             \bar{x}_j & \hbox{ otherwise}
             \end{array}
  \right.\]
Once again, following a dual reasoning to the proof of Theorem~\ref{th:sys_maximal} and considering Statement (b) in Lemma~\ref{lemma:sys_solutionset}, we conclude that $f$ is a bijection. Hence, the number of minimal solutions of Equation~\eqref{sys:bipolar_general} coincides with the number of maximal elements of $S^-$.
\qed\end{proof}

From Theorem~\ref{th:sys_minimal}, we can deduce the analytical description of the least solution or minimal solutions of a bipolar max-product FRE. The next corollary formalizes this result.

\begin{corollary}\label{cor:sys_minimal}
{Let Equation~\eqref{sys:bipolar_general} be a solvable bipolar max-product FRE, $(\bar{x}_1,\bar{y}_1,\dots,\bar{x}_m,\bar{y}_m)\in[0,1]^{2m}$ be the greatest solution of its corresponding max-product FRE~\eqref{sys:FRE_general}, $S$  be the set of feasible pairs with respect to Equation~\eqref{sys:FRE_general} and $S^-=\{J^-\mid(J^+,J^-)\in S\}$. Consider the mapping $f\colon S^-\to [0,1]^m$ which associates each $J^-\in S^-$ with the tuple $(x_1,\dots,x_m)$ defined, for each $j\in\{1,\dots,m\}$, as:
\[x_j= \left\{\begin{array}{ll}
             1-\bar{y}_j & \hbox{ if}\quad j\in J^- \\
             \bar{x}_j & \hbox{ otherwise}
             \end{array}
  \right.\]
The following statements hold:
  \begin{enumerate}
  \item[(a)] If the set $S^-$ has a greatest element $J^-$, then $f(J^-)$ is the least solution of Equation~\eqref{sys:bipolar_general}.
  \item[(b)] Let $M^-$ be the set of maximal elements of $S^-$. Then, the set of minimal solutions of Equation~\eqref{sys:bipolar_general} is given by:
  \[\{f(J^-)\mid J^-\in M^-)\}\]
  \end{enumerate}
  }
\end{corollary}

Finally, we will introduce an illustrative example in order to clarify the results provided in this section.

\begin{example}\label{ex:bipolarFRE}
The following bipolar max-product FRE composed of three equations and three unknown variables will be considered:
{\small\begin{equation}\label{eq:bipolarFRE_ex}
\begin{array}{ll}
(0.1*x_1) \vee (0.3*(1-x_1)) \vee (0.25*x_2) \vee (0.3*(1-x_2)) \vee (0.4*x_3) \vee (0.4*(1-x_3))&=0.2\\
(0.1*x_1) \vee (0.8*(1-x_1)) \vee (0.3*x_2) \vee (0.6*(1-x_2)) \vee (0.3*x_3) \vee (0.5*(1-x_3))&=0.4\\
(0.3*x_1) \vee (0.8*(1-x_1)) \vee (0.4*x_2) \vee (0.5*(1-x_2)) \vee (0.8*x_3) \vee (0.8*(1-x_3))&=0.4
\end{array}
\end{equation}}
It is easy to see that the corresponding max-product fuzzy equation of Equation~\eqref{eq:bipolarFRE_ex} is given by {Equation~\eqref{eq:FRE_ex} in Example~\ref{ex:feasiblepair}}.

In the following, we will check whether the conditions required in Theorem~\ref{th:sys_general_mxn} are satisfied. Clearly, the inequality $1\leq \bar{x}_i+\bar{y}_i$ holds, for each $i\in\{1,2,3\}$. {Furthermore, as shown in Example~\ref{ex:feasiblepair}, the pair $(\{2\},\{1,3\})$ forms a feasible pair with respect to Equation~\eqref{eq:FRE_ex}.}

As a result, the hypothesis required in Theorem~\ref{th:sys_general_mxn} are verified and therefore, we can assert that there exists at least a solution of Equation~\eqref{eq:bipolarFRE_ex}. For instance, we can deduce that the tuple $(0.5,0.8,0.5)$ is a solution of such system.

Now, making the corresponding calculations, we obtain that the sets $S^+$ and $S^-$ are given by:
\begin{equation*}
\begin{array}{ll}
S^+=&\Big\{\varnothing,\{1\},\{2\},\{3\},\{1,3\},\{2,3\}\Big\}
\vspace{0.1cm}\\
S^-=&\Big\{\{1\},\{2\},\{1,2\},\{1,3\},\{2,3\},\{1,2,3\}\Big\}
\end{array}
\end{equation*}
It is important to highlight that there is no $A\in S^+$ such that $\{1,2\}\subseteq A$. The reason is that, in such case, the second equation is not satisfied. Furthermore, this fact implies that $\{1,2,3\}\notin S^+$ and consequently $\varnothing\notin S^-$.

On the one hand, notice that $S^+$ has two maximal elements, $\{1,3\}$ and $\{2,3\}$. As a consequence, according to Theorem~\ref{th:sys_maximal}, Equation~\eqref{eq:bipolarFRE_ex} has two maximal solutions. In particular, Corollary~\ref{cor:sys_maximal} allows us to ensure that $(\bar{x}_1,1-\bar{y}_2,\bar{x}_3)=(1,0.\widehat{3},0.5)$ and $(1-\bar{y}_1,\bar{x}_2,\bar{x}_3)=(0.5,0.8,0.5)$ are the two maximal solutions of Equation~\eqref{eq:bipolarFRE_ex}. 

On the other hand, as $\{1,2,3\}$ is the greatest element of $S^-$, applying Theorem~\ref{th:sys_minimal} we deduce that Equation~\eqref{eq:bipolarFRE_ex} has a least solution. Specifically, from Corollary~\ref{cor:sys_minimal}, we conclude that the tuple $(1-\bar{y}_1,1-\bar{y}_2,1-\bar{y}_3)=(0.5,0.\widehat{3},0.5)$ is the least solution of Equation~\eqref{eq:bipolarFRE_ex}.
\qed
\end{example}

\section{Conclusions and future work}\label{sec:conclusion} 
We have introduced {a general study} on the solvability of bipolar max-product fuzzy relation equations with the standard negation. First of all, we have analysed the particular case of bipolar max-product fuzzy equations. As well as characterizing the solvability of these equations, we have included different properties related to the existence of a greatest/least solution or maximal/minimal solutions. 

In the sequel, the necessary and sufficient conditions, which indicate when a bipolar max-product fuzzy relation equation with the standard negation is solvable, have been presented. Finally, we have studied the maximal and minimal solutions of a solvable bipolar max-product fuzzy relation equation with the standard negation. This study allows to solve equations with two of the most interesting and useful operators considered in real cases, the product t-norm (related to probabilities properties) and the standard negation (the most used negation in the applications). {To complete the presented investigation on bipolar max-product FREs with the standard negation, providing the set of all solutions will suppose one of the main challenges to be explored in the future. Furthermore, a formal study on the modification of the covering problem for bipolar FREs, as well as on the concept of irredundant covering for bipolar FREs, will be carried out in the future.}

{In order to provide an automatic mechanism for the solvability of a bipolar max-product fuzzy relation equation with the standard negation, an efficient algorithm for the existence of feasible pairs will be researched in the future. For this purpose, existing algorithms for FREs in the literature will be taken into account as reference~\cite{BOURKE1998,Loetamonphong:99,Luoh2002,Pedrycz1984,Peeva2013,Peeva2007,Zahariev:2010}. In particular, due to the connection between feasible pairs and coverings, those algorithms concerning the covering problem~\cite{Lin2011,Markovskii,Shieh2013,Yeh200823} will be instrumental in this task.} In addition to applying the obtained results to practical examples, we are interested in studying bipolar fuzzy relation equations based on more general compositions such as max-t-norm, max-uninorm and max-unorm. The use of arbitrary negations in bipolar fuzzy relation equations will also be an important issue to deal with in our future work.

\end{document}